\documentclass{amsart}

\newtheorem{theorem}{Theorem}[section]
\newtheorem{lemma}[theorem]{Lemma}

\theoremstyle{definition}
\newtheorem{definition}[theorem]{Definition}

\theoremstyle{remark}
\newtheorem{remark}[theorem]{Remark}

\numberwithin{equation}{section}

\begin{document}

\title{Deformations of $2k$-Einstein structures}

%    Information for first author
\author{Levi Lopes de Lima}
%    Address of record for the research reported here
\address{Federal University of Cear\'a,
Department of Mathematics, Campus do Pici, R. Humberto Monte, s/n, 60455-760,
Fortaleza/CE, Brazil. }
%    Current address
%\curraddr{Department of Mathematics and Statistics,
%Case Western Reserve University, Cleveland, Ohio 43403}
\email{levi@mat.ufc.br}
%    \thanks will become a 1st page footnote.
\thanks{The first author was supported in part by a CNPq  grant.}

%    Information for second author
\author{Newton Lu\'{\i}s Santos}
\address{Federal University of
Piau\'{\i},
 Department of Mathematics, Campus Petronio Portela, Ininga, 64049-550 Teresina/PI, Brazil}
\email{newtonls@ufpi.br}
\thanks{The second author was supported in part by the  CNPq postdoctoral grant 155723/2006-5.}

%    General info

%\date{January 1, 2001 and, in revised form, June 22, 2001.}


\begin{abstract}
It is shown that the space of infinitesimal deformations of $2k$-Einstein structures is finite dimensional at compact non-flat
space forms. Moreover, spherical space forms are shown to be rigid in the sense that they are isolated in the corresponding moduli space.
\end{abstract}

\maketitle

\section{Introduction and statement of results}

Let $X$ be a
compact connected manifold of dimension $n\geq 3$. We  denote by $\Omega^p(X)$ the
space of exterior $p$-forms on $X$. The space of double forms
of bi-degree $(p,q)$ is
$\Omega^{p,q}(X)=\Omega^p(X)\otimes\Omega^q(X)$, the tensor product
taken over the ring of smooth functions on $X$. Clearly,
$\Omega(X)=\oplus_{p,q\geq 0}\Omega^{p,q}(X)$ is a bigraded associative algebra.
For example, any bilinear form on tangent vectors is a
$(1,1)$-form. In particular, if ${\mathcal M}(X)$ is the space of Riemannian metrics on $X$ then any $g\in{\mathcal M}(X)$ is a
$(1,1)$-form. Also, the curvature tensor $R_g$ of $g$ can be viewed as a
$(2,2)$-form. In fact, if we define ${\mathcal C}^p(X)\subset \Omega^{p,p}(X)$ to be the space of $(p,p)$-forms satisfying the symmetry condition
$$
\omega(x_1\wedge\ldots\wedge x_{p}\otimes y_1\wedge\ldots\wedge
y_{p})=\omega(y_1\wedge\ldots\wedge y_{p}\otimes x_1\wedge\ldots\wedge
x_{p}),
$$
then any bilinear symmetric form on tangents vectors ($g$ in particular) belongs   to ${\mathcal C}^1(X)$, and moreover $R_g\in{\mathcal C}^2(X)$ by its symmetries. An account of the formalism of double forms used in this paper can be found in \cite{L1}, \cite{L2}.

Let us consider a Riemannian manifold $(X,g)$, with $X$ as above. Then multiplication
by the metric defines a map $g:\Omega^{p-1,q-1}(X)\to\Omega^{p,q}(X)$.
Also, a contraction operator $c_g:\Omega^{p,q}(X)\to\Omega^{p-1,q-1}(X)$ is defined by
$$
c_g\omega(x_1\wedge\ldots\wedge x_{p-1}\otimes y_1\wedge\ldots\wedge
y_{p-1})=\sum_i \omega(e_i\wedge x_1\wedge\ldots\wedge
x_{p-1}\otimes e_i\wedge y_1\wedge\ldots\wedge y_{p-1}),
$$
where $\{e_i\}$ is a local orthonormal tangent frame. Clearly, $g$ and $c_g$
are tensorial and it can be shown that, considered as operators on $\Lambda^{*,*}_x$, $x\in X$, they are adjoint to each
other  with respect to the natural inner product.

\begin{remark}\label{const}
In the language of double forms, that a Riemannian manifold $(X,g)$ has
constant sectional curvature $\mu$ can be characterized  by the identity
$R_g=\frac{\mu}{2}g^2$.
\end{remark}

Contraction can be used to rewrite the Ricci tensor ${\mathcal R}_g$ and the
scalar curvature ${\mathcal S}_g$ of $(X,g)$ as ${\mathcal R}_g=c_gR_g$ and ${\mathcal S}_g=c^2_gR_g$. More
generally, we set for any $1\leq k\leq n/2$,
\begin{equation}\label{def}
{\mathcal R}^{(2k)}_g=c^{2k-1}_gR^k_g,\quad {\mathcal S}^{(2k)}_g=\frac{1}{(2k)!}c^{2k}_gR^k_g.
\end{equation}
These are the $2k$-Ricci tensor and the $2k$-curvature, respectively. We also consider the $2k$-Einstein-Lovelock tensor \cite{L2} \cite{Lo}
$$
{\mathcal T}^{(2k)}_g=\frac{{\mathcal R}^{(2k)}_g}{(2k-1)!}-{\mathcal S}^{(2k)}_g g.
$$
For $k=1$ we recover the standard notions of Ricci tensor, (half) the scalar curvature and the Einstein tensor. We remark that, as shown in \cite{Lo}, the Einstein-Lovelock tensors span the subspace of ${\mathcal C}^1$ formed by divergence free elements which depend on the derivatives of $g$ up to second order. Also, these geometric invariants play a fundamental role in stringy gravity \cite{C-B}.

Following \cite{L2} we say that $(X,g)$ is $2k$-Einstein if there exists a smooth function $\lambda$ on $X$ such that
$$
{\mathcal R}_g^{(2k)}=\lambda g.
$$
Thus, $2$-Einstein means precisely that $(X,g)$ is Einstein in the usual sense. It is  shown in \cite{L2} that $2k$-Einstein metrics are critical for the Einstein-Hilbert-Lovelock functional
$$
{\mathcal A}^{(2k)}(g)=\int_X {\mathcal S}^{(2k)}_g\nu_g,
$$
restricted to the space ${\mathcal M}_1(X)$ of metrics of unit volume. Here, $\nu_g$ is the volume element of $g$.  Examples of $2k$-Einstein manifolds include spaces forms and isotropically irreducible homogeneous Riemannian manifolds. Also, if $2k=n$ then any metric is $2k$-Einstein, since ${\mathcal S}^{(n)}_g$ is, up to a constant, the Gauss-Bonnet integrand. Thus we may assume $n>2k$ in what follows.
Under this assumption a standard argument using the fact that  ${\mathcal T}^{(2k)}_g$ is divergence free implies that  $\lambda$ is actually a constant; in fact, since ${\mathsf {tr}}_g {\mathcal R}^{(2k)}_g=\langle c^{2k-1}_gR^k_g,g\rangle=c^{2k}_gR^k_g=(2k)!{\mathcal S}^{(2k)}_g$, we have
\begin{equation}\label{const}
\lambda=\frac{(2k)!}{n}{\mathcal S}^{(2k)}_g,
\end{equation}
so that ${\mathcal S}^{(2k)}_g$ is constant as well.

At this point we follow \cite{B} and appeal to a theorem of Moser \cite{M}: if $g,g'\in{\mathcal M}_1(X)$ there exists a diffeomorphism $\varphi:X\to X$ such that $\nu_{g'}=\varphi^*\nu_g$. Thus, in  order to understand
the structure of the moduli space ${\mathcal E}_{2k}(X)$ of $2k$-Einstein metrics on $X$, after moding out by the standard action of the diffeomorphism group $\mathfrak{D}(X)$ of $X$ on ${\mathcal M}_1(X)$, we may restrict ourselves to the space ${\mathcal N}_{\nu}=\{g\in{\mathcal M}(X);\nu_g=\nu\}$ of metrics with a fixed volume element $\nu$. If $g$ is $2k$-Einstein, the corresponding element in ${\mathcal E}_{2k}(X)$ will be represented by $[g]$, and will be referred to as a $2k$-Einstein structure. We note that $2k$-Einstein metrics are critical points for ${\mathcal A}^{(2k)}$ restricted to ${\mathcal N}_{\nu}$.

A first step toward understanding the structure of ${\mathcal E}_{2k}(X)$ is to estimate the dimension of the space of infinitesimal deformations of $2k$-Einstein structures at a given $[g]\in{\mathcal E}_{2k}(X)$. To motivate this approach, consider
first order jets $h=dg_t/dt|_{t=0}$ of one-parameter families (i.e. {\em deformations}) $g_t\in{\mathcal N}_{\nu_g}$ of $2k$-Einstein metrics with $g_0=g$. Naturally, here we should discard {\em trivial} deformations of the type  $g_t=\varphi^*_tg$ for a one parameter family of diffeomorphisms with $\varphi_0={\mathsf {Id}}_X$. It is known \cite{B} that these jets are characterized by $h\in {\mathsf {Im}}\,\, \delta^*_g$, where $\delta^*_g:\Omega^1(X)\to {\mathcal C}^1(X)$ is the $L^2$ adjoint of $\delta_g:{\mathcal C}^{1}(X)\to\Omega^1(X)$, the standard divergence operator acting on symmetric $2$-tensors. Thus, by results in \cite{BE}, essential deformations necessarily satisfy $\delta_gh=0$, which means that $h$ is orthogonal to the tangent space at $g$ of the orbit of $g$ with respect to the action of $\mathfrak{D}(X)$ on ${\mathcal M}_1(X)$. Also, since in general there holds
\begin{equation}\label{mean}
\dot\nu_gh\doteq\frac{d}{dt}\nu_{g_t}|_{t=0}=\frac{1}{2}{\mathsf {tr}_gh}\,\nu_g,
\end{equation}
we must add the condition ${\mathsf {tr}}_gh=0$, which  is precisely the infinitesimal way of saying that $g_t\in{\mathcal N}_{\nu_g}$. In what follows we shall retain the dot notation for the linearization of other invariants of $g$.

\begin{remark}\label{expdiv}
Our sign convention for $\delta_g$ is
$$
(\delta_g\omega)(x) = - \sum_i(\nabla_{e_i}\omega)(e_i,x),
$$
where $\nabla$ is the covariant derivative on ${\mathcal C}^1(X)$ induced by the Levi-Civita connection of $g$, also denoted by $\nabla$.
\end{remark}

Notice also that if $g_t$ is a deformation
of $2k$-Einstein metrics with $g_0=g$ then we have ${\mathcal R}^{(2k)}_{g_t}=\lambda_{g_t}g_t$, which yields, after (\ref{const}),
$$
\dot{\mathcal R}^{(2k)}_gh=\frac{(2k)!}{n}\dot{\mathcal S}^{(2k)}_ghg+\lambda_gh.
$$
Now, each ${\mathcal S}^{(2k)}_{g_t}$ is constant over $X$ so that ${\mathcal S}^{(2k)}_{g_t}={\mathcal A}^{(2k)}(g_t)$, and since $g_t$ is critical for ${\mathcal A}^{(2k)}$ restricted to ${\mathcal N}_{\nu_g}$, the
above discussion justifies the following definition.

\begin{definition}
If $(X,g)$ is $2k$-Einstein, ${\mathcal R}^{(2k)}_g=\lambda_g g$, the space of infinitesimal $2k$-Einstein deformation at $[g]$, denoted $\mathfrak{M}_{[g]}^{(2k)}$, is the vector space of all $h\in{\mathcal C}^1(X)$ such that
\begin{equation}\label{moduli}
\dot{\mathcal R}_g^{(2k)}h=\lambda_g h,
\end{equation}
and
\begin{equation}\label{moduli2}
\quad \delta_gh=0,\quad {\mathsf {tr}}_gh=0.
\end{equation}
\end{definition}

Thus, if we set ${\mathcal I}_g=\delta_g^{-1}(0)\cap{\mathsf {tr}}^{-1}_{g}(0)$ and ${\mathcal L}_g^{(2k)}=\dot{\mathcal R}_g^{(2k)}-\lambda_g $ then $\mathfrak{M}_{[g]}^{(2k)}=\ker {\mathcal L}_g^{(2k)}|_{{\mathcal I}_g}$. In particular, since $X$ is compact, $\mathfrak{M}_{[g]}^{(2k)}$ is finite dimensional if ${\mathcal L}_g^{(2k)}|_{{\mathcal I}_g}$
is elliptic.

A theorem of Ebin \cite{E} shows that ${\mathcal I}_g$ can be \lq exponentiated\rq\, to yield a local slice $\mathfrak{S}_g$ for the $\mathfrak{D}(X)$-action.

\begin{definition}\label{premod}
The space of all $2k$-Einstein metrics in $\mathfrak{S}_g$ is called the {\em pre-moduli space around} $g$, denoted ${\mathcal E}_{2k}(X)_g$.
\end{definition}

The moduli space itself ${\mathcal E}_{2k}(X)$ is locally obtained after dividing ${\mathcal E}_{2k}(X)_g$ by the isometry group of $g$. In what follows, however, we shall completely ignore this issue and  work with ${\mathcal E}_{2k}(X)_g$ directly.

It has been shown in \cite{BE} that ${\mathcal L}_g^{(2)}|_{{\mathcal I}_g}$ is elliptic for any $(X,g)$ Einstein, thus showing affirmatively the finite dimensionality of $\mathfrak{M}_{[g]}^{(2k)}$ in the case $k=1$. If $k\geq 2$, however,
this is not generally true as $\mathfrak{M}_{[g]}^{(2k)}$ may be infinite dimensional for certain choices of $(X,g)$, reflecting the fact that ${\mathcal L}_g^{(2k)}|_{{\mathcal I}_g}$ might be of mixed type. For example, consider the Riemennian product $X=M^r\times T^m$, where $M$ is an arbitrary Riemannian manifold and $T^m$ is a flat torus. If  $2k>r$ then $X$ is $2k$-Einstein, irrespective of the metric on $M$, thus showing that $\dim\mathfrak{M}_{[g]}^{(2k)}=+\infty$ in this case.
In view of this example, it is an interesting  to find $2k$-Einstein structures $(X,g)$ for which $\mathfrak{M}_{[g]}^{(2k)}$ is finite dimensional. In this respect, we prove here:

\begin{theorem}\label{main}
If $(X,g)$ is a compact non-flat space form then $\mathfrak{M}_{[g]}^{(2k)}$ is finite dimensional. Moreover, if $(X,g)$ is spherical then $\mathfrak{M}_{[g]}^{(2k)}$ is trivial.
\end{theorem}

Thus, adapting the terminology in \cite{K1}, we see that spherical space forms are infinitesimally $2k$-non-deformable.

Again in accordance with \cite{K1} we say that $[g]$ is $2k$-{\em non-deformable} if any deformation $g_t$ of $g$ by $2k$-Einstein metrics is trivial.
Using Theorem \ref{main} and a result due to Koiso \cite{K1} we can check without effort  that spherical space forms are indeed $2k$-non-deformable. However we can adapt an argument in \cite{B} to show a much stronger result, namely, that spherical space forms are rigid.

\begin{theorem}\label{main2}
If $(X,g)$ is a compact spherical space form then $[g]$ is an isolated point in ${\mathcal E}_{2k}(X)_g$.
\end{theorem}

In particular, the only way of locally deforming the standard metric on the unit sphere by $2k$-Einstein metrics is by dilations. This should be compared to an old result by Berger \cite{Be}, where rigidity results have been proved for sufficiently positively pinched Einstein structures.

In fact, Theorem \ref{main2}  will be deduced here from a more detailed structure result for the moduli space near space forms.

\begin{theorem}\label{main3}
Near to a compact non-flat space form $(X,g)$, ${\mathcal E}_{2k}(X)_g$ has the structure of a real analytic subset which in turn
is contained in a finite dimensional analytic manifold whose tangent space at $[g]$ is precisely $\mathfrak{M}_{[g]}^{(2k)}$.
\end{theorem}

\section{Some preliminary facts}

The proof of Theorems \ref{main} and  \ref{main2} rely on several facts established in \cite{L1} and \cite{L2} that we recall in this section. We start with a  commutation formula for powers of $g$
and $c_g$. It is easy to check that
\begin{equation}\label{kulk}
c_gg\omega=gc_g\omega+(n-p-q)\omega,\qquad \omega\in\Omega^{p,q}(X).
\end{equation}
More generally, it is proved in \cite{L1} that, for $\omega\in\Omega^{p,q}(X)$,
\begin{equation}\label{comut}
\frac{1}{m!}c^l_gg^m\omega=\frac{1}{m!}g^mc^l_g\omega+\sum_{r=1}^{\min \{l,m\}} C^l_r\prod_{i=0}^{r-1}(n-p-q+l-m-i)\frac{g^{m-r}}{(m-r)!}c^{l-r}_g\omega,
\end{equation}
where $C^l_r$ is the standard binomial coefficient.
This identity will reveal itself to be  extremely important in the sequel.

Also, for
$h\in {\mathcal C}^1(X)$ it is defined in \cite{L2}, the linear map $F_h:{\mathcal C}^p(X)\to {\mathcal C}^p(X)$,
\begin{eqnarray*}
&&F_h\omega(e_{i_1}\wedge \ldots \wedge e_{i_p},e_{j_1}\wedge \ldots
\wedge
e_{j_p})=\\
&&\qquad
\Big(\sum_{k=1}^ph(e_{i_k},e_{i_k})+\sum_{k=1}^ph(e_{j_k},e_{j_k})
\Big)\omega(e_{i_1}\wedge \ldots \wedge e_{i_p},e_{j_1}\wedge \ldots
\wedge e_{j_p}).
\end{eqnarray*}
It turns out that $F_h$ has many interesting properties (it is self-adjoint, acts as a
derivation, etc.) but here we will only need the following ones.

\begin{itemize}
\item $F_h(g^p)=2pg^{p-1}h$ for $p\ge 1$;
\item If, as before, $R_g$ is the curvature tensor of $g$,
\begin{eqnarray}\label{rrr}
F_h(R_g)(x\wedge y,z\wedge
u)& = & h(R_g(x,y)z,u)-h(R_g(x,y)u,z)+\nonumber\\
&  & \quad +h(R_g(z,u)x,y)-h(R_g(z,u)y,x),
\end{eqnarray}
where in the right-hand side $R_g$ is viewed as a
$(3,1)$-tensor.
\end{itemize}

Also, consider the second order differential operator ${\mathcal D}^2:{\mathcal C}^1(X)
\to {\mathcal C}^2(X)$ given by:
\begin{eqnarray*}
 {\mathcal D}^2 h (x_1\wedge x_2,y_1\wedge y_2) & = & \nabla^2_{y_1,x_1}h(x_2,
y_2)+ \nabla^2_{x_1,y_1}h(x_2, y_2) \\
& & +\nabla^2_{y_2,x_2}h(x_1,
y_1)
+\nabla^2_{x_2,y_2}h (x_1, y_1)\\
& &-\nabla^2_{y_1,x_2}h (x_1,
y_2)-
\nabla^2_{x_2,y_1}h (x_1, y_2)\\
& & -\nabla^2_{y_2,x_1}h(x_2,
y_1)
-\nabla^2_{x_1,y_2}h (x_2, y_1),
\end{eqnarray*}
where
$$
\nabla^2_{x,y}h(w,z)=\nabla_x(\nabla_yh)(w,z)-(\nabla_{\nabla_xy}h)(w,z).
$$
Then it is shown in
\cite{L2} that the linearization of the tensor curvature, projected on ${\mathcal C}^2(X)$, is given
by
\begin{equation}\label{varl}
\dot R_g h=-\frac14 {\mathcal D}^2 h+\frac14 F_h(R_g).
\end{equation}

Later on we will show that, on a $2k$-Einstein manifold with constant curvature, the second order part of ${\mathcal L}_g^{(2k)}$ is completely determined by the first and second contractions of $\dot R_g$. Thus, in view of (\ref{varl}), we  need to determine the first and second contractions of ${\mathcal D}^2$ and $F_h(R_g)$. The result, for a general metric $g$, appears in the following lemma, whose formulation requires some more notation.

\begin{definition}\label{defn}
If $\xi,\eta\in{\mathcal C}^1(X)$  we set
\begin{equation}\label{prod}
(\xi\circ\eta)(v,w)=\sum_{i=1}^n \xi(v,e_i)\eta(e_i,w).
\end{equation}
Moreover, ${\mathfrak{R}}:{\mathcal C}^1(X)\to {\mathcal C}^1(X)$ is defined by
$$
({\mathfrak{R}} h)(v,w)=\sum_{i=1}^n h(R(v,e_i)w,e_i).
$$
\end{definition}

\begin{lemma}\label{aux}
For a general metric $g$ on $X$ the following identities hold:
\begin{enumerate}
 \item if $\nabla^*\nabla$ is the Bochner Laplacian on ${\mathcal C}^1(X)$,
  \begin{eqnarray}\label{cR}
c_g{\mathcal D}^2 h & = & -2(\nabla^*\nabla h)+2
\nabla^gd\mathsf {tr}_g h+4\delta^{*}_g\delta_g h-\nonumber\\
& & \qquad \qquad-({\mathcal R}_g\circ h+h\circ {\mathcal R}_g) +2{\mathfrak R} h.
\end{eqnarray}
\item if $\Delta_g$ is the Laplacian of $g$,
\begin{equation}\label{c2R}
c^2_g({\mathcal D}^2 h) = -4\Delta_g \mathsf {tr}_gh-4\delta_g \delta_g h.
\end{equation}
\item
$c_gF_h(R_g)={\mathcal R}_g\circ h+h\circ {\mathcal R}_g+ 2{\mathfrak R}h.$
\item
$c^2_gF_h(R_g)=4\langle{\mathcal R}_g,h\rangle.$
\end{enumerate}
\end{lemma}

The proofs of items 3. and 4. are rather easy computations. In fact,
using (\ref{rrr}) and Definition \ref{defn}, we compute:
\begin{eqnarray*}
c_gF_h(R_g)(x,y)&=&\sum_i F_h(R_g)(e_i\wedge x,e_i\wedge
y)\\
&=& \sum_i\Big(h(R_g(e_i,x)e_i,y)-h(R_g(e_i,x)y,e_i)+\\
 & & \qquad + h(R_g(e_i,y)e_i,x)-h(R_g(e_i,y)x,e_i)\Big)\\
&=& ({\mathcal R}_g\circ h+h\circ {\mathcal R}_g)(x,y)+ 2({\mathfrak R}h) (x,y).
\end{eqnarray*}
Moreover,
$$
c^2_gF_h(R_g)=c_g({\mathcal R}_g\circ h+h\circ {\mathcal R}_g+ 2{\mathfrak{R}h})=4\langle{\mathcal R}_g,h\rangle.
$$

The proofs of the remaining items are a bit more involved. In fact, we have
\begin{eqnarray*}
c_g({\mathcal D}^2 h)(e_j,e_k) & = & \sum_i{\mathcal D}^2 h( e_i\wedge e_j,e_i\wedge e_k)\\
& = &\sum_i\Big[
2\nabla^2_{e_i,e_i}h(e_j, e_k)+\left( \nabla^2_{e_k,e_j}h(e_i, e_i)
+\nabla^2_{e_j,e_k}h(e_i, e_i) \right)-\\
& & \quad -\left( \nabla^2_{e_k,e_i}h(e_j, e_i) + \nabla^2_{e_j,e_i}h (e_i,
e_k) \right)-\\
& & \qquad -\left( \nabla^2_{e_i,e_j}h(e_i, e_k)
+\nabla^2_{e_i,e_k}h(e_j, e_i)\right)\Big].
\end{eqnarray*}
We can use  Ricci identities to rewrite the last term above as
\begin{eqnarray*}
\nabla^2_{e_i,e_j}h(e_i, e_k)
+\nabla^2_{e_i,e_k}h(e_i,e_j) & = & \nabla^2_{e_j,e_i}h(e_i, e_k)
+\nabla^2_{e_k,e_i}h(e_i,e_j)+\\
& &  +h(R(e_i,e_j)e_i,e_k)
+h(e_i,R(e_i,e_j)e_k)+ \\
& & +h(R(e_i,e_k)e_i,e_j) +h(e_i,R(e_i,e_k)e_j),
\end{eqnarray*}
so that
\begin{eqnarray*}
c_g({\mathcal D}^2 h)(e_j,e_k) & = &\sum_i\Big[ 2\left(\nabla^2_{e_i,e_i}h\right)(e_j, e_k)+\left( \nabla^2_{e_k,e_j}h +\nabla^2_{e_j,e_k}h \right)(e_i, e_i)-\\
& & -2\Big( \nabla^2_{e_k,e_i}h(e_i,e_j) + \nabla^2_{e_j,e_i}h (e_i, e_k) \Big)-\\
& & -\Big( h(R(e_i,e_j)e_i,e_k) +h(e_i,R(e_i,e_j)e_k)+\\
& &+h(R(e_i,e_k)e_i,e_j) +h(e_i,R(e_i,e_k)e_j) \Big)\Big]\\
&= &-2(\nabla^*\nabla h)(e_j, e_k)+\\
& &+\Big( \nabla d\mathsf {tr}_gh(e_k,e_j) +\nabla d\mathsf {tr}_gh(e_j,e_k)\Big)-\\
&& +4\delta^{*}_g\delta_g h(e_k,e_j)-({\mathcal R}_g\circ h)(e_j,e_k)-\\
& & -(h\circ {\mathcal R}_g)(e_j,e_k) +2({\mathfrak R} h)(e_j,e_k),
\end{eqnarray*}
and the first item follows. Finally,
using the symmetries due to the fact that ${\mathcal D}^2 h\in {\mathcal C}^2(X)$ we
obtain
\begin{eqnarray*}
c^2_g{\mathcal D}^2 h &  = & \sum_{ij}{\mathcal D}^2 h(e_i\wedge e_j,e_i\wedge e_j)\\
& = & \sum_{ij}\left(
4(\nabla^2_{e_i,e_i}h)(e_j, e_j)-4(\nabla^2_{e_i,e_j}h )(e_i, e_j)\right).
\end{eqnarray*}
Now notice  that
\begin{equation}
\sum_{ij}\nabla^2_{e_i,e_i}h(e_j, e_j)=-\Delta_g \mathsf {tr}_g(h), \qquad
\sum_{ij} \nabla^2_{e_i,e_j}h (e_i,e_j)=\delta_g (\delta_g h),
\end{equation}
and the second item follows. Lemma \ref{aux} is proved.

\begin{remark}\label{newr}
Under the constant curvature condition $R_g=\frac{\mu}{2}g^2$ (see Remark \ref{const}) we have ${\mathcal R}_g=(n-1)\mu g$, so that
\begin{equation}
{\mathcal R}_g\circ h+h\circ {\mathcal R}_g=(n-1)\mu(g\circ h+h\circ g)=2(n-1)\mu
h \label{dd}
\end{equation}
and
\begin{equation}
{\mathfrak R} h= \mu ({\rm tr}_gh \,g - h) \label{ee}.
\end{equation}
\end{remark}

\section{The proof of Theorem \ref{main}}

As already observed in Remark \ref{const}, if $(X,g)$ has constant sectional curvature $\mu$ then its curvature tensor can be expressed as $R_g=\frac{\mu}2g^2$, so that  its $2k$-Ricci tensor can be computed using (\ref{comut}) as follows:
\begin{eqnarray*}
{\mathcal R}^{(2k)}_g&=&c^{2k-1}_g\left(\frac{\mu^k}{2^k}g^{2k}\right)\\
& = & \frac{\mu^k}{2^k}c^{2k-1}_gg^{2k-1}g\\
&=&\frac{\mu^k}{2^k}(2k-1)!\sum_{r=1}^{2k-1}C^{2k-1}_r
\prod_{i=0}^{r-1}(n-2-i)\frac{g^{2k-1-r}}{(2k-1-r)!}c^{2k-1-r}_gg\\
&=&\frac{\mu^k}{2^k}(2k)!\prod_{i=0}^{2k-3}(n-2-i)(n-1)g\\
 & = & \frac{(n-1)!(2k)!}{(n-2k)!2^k}\mu^kg.
\end{eqnarray*}
Thus, $g$ is a $2k$-Einstein metric,
$
{\mathcal R}^{(2k)}_g=\lambda_k g$,
with
\begin{equation}\label{eins}
\lambda_k=(n-1)(n-2)C_{n,k}\mu^k,\qquad C_{n,k}=\frac{(2k)!(n-3)!}{2^k(n-2k)!}.
\end{equation}

From (\ref{def}) we have in general
\begin{equation}\label{sy}
\dot{\mathcal R}^{(2k)}_gh = (2k-1)(\dot c_gh)c^{2k-2}_gR^k_g +
kc^{2k-1}_gR^{k-1}_g\dot R_gh,
\end{equation}
which specializes to
\begin{eqnarray}
\dot{\mathcal R}^{(2k)}_gh & = & \frac{(2k-1)\mu^k}{2^k}(\dot c_gh)c^{2k-2}g^{2k}
+\frac{k\mu^{k-1}}{2^{k-1}}c^{2k-1}_gg^{2k-2}\dot R_gh\nonumber\\
& \doteq & A_gh +B_gh\label{var1}
\end{eqnarray}
in the constant curvature case.
We now identify the terms $A_gh$ and $B_gh$ in this expression.

We start with $A_gh$, which should be a zero order term in $h$ (no derivatives).
To begin with
observe that (\ref{comut}) implies
\begin{eqnarray*}
c^{2k-2}_g\left(\frac{g^{2k-1}}{(2k-1)!}g\right) & = &\sum_{r=1}^{2k-2}C^{2k-2}_r
\prod_{i=0}^{r-1}(n-3-i)\frac{g^{2k-1-r}}{(2k-1-r)!}c^{2k-2-r}_gg\\
& = & (2k-2) \prod_{i=0}^{2k-4}(n-3-i)\frac{g^2}2c_gg+
\prod_{i=0}^{2k-3}(n-3-i)g^2\\
& = &
\prod_{i=0}^{2k-4}(n-3-i)\Big(n(k-1)+(n-2k)\Big)g^2\\
& = &\frac{(n-2)!k}{(n-2k)!}~g^2,
\end{eqnarray*}
that is
\begin{equation}
c^{2k-2}_gg^{2k}= \frac{(n-2)!(2k)!}{2(n-2k)!}~g^2,
\end{equation}
which gives
\begin{eqnarray*}
A_gh=
\frac{(2k-1)(n-2)!(2k)!\mu^k}{(n-2k)!2^{k+1}}(\dot c_g h)g^2.
\end{eqnarray*}
Now if we linearize the identity  $c_gg^2=2(n-1)g$ in the direction of  $h$ we obtain, with the help of
(\ref{kulk}),
\begin{eqnarray*}
2(n-1)h&=&(\dot c_gh)g^2+2c_ggh\\
&=&(\dot c_gh)g^2+2{\mathsf {tr}}_gh\,g+2(n-2)h,
\end{eqnarray*}
so that
$$
(\dot c_gh)g^2=2\left(h-{\mathsf {tr}}_ghg\right),
$$
implying
\begin{eqnarray}
A_gh=
(2k-1)(n-2)C_{n,k}\mu^k(h-{\mathsf {tr}}_ghg). \label{zz}
\end{eqnarray}

In order to compute $B_gh$ we make use of
(\ref{comut}) with
$k>1$, an assumption implying in particular that $c^{2k-1}_g\dot R_gh=0$ since $\dot R_gh\in {\mathcal C}^2(X$). Hence,
\begin{eqnarray*}
c^{2k-1}_g\left(\frac{g^{2k-2}}{(2k-2)!}\dot R_gh\right) & = & \sum_{r=1}^{2k-2}C^{2k-1}_r
\prod_{i=0}^{r-1}(n-3-i)\frac{g^{2k-2-r}}{(2k-2-r)!}c^{2k-1-r}_g\dot R_gh\\
& = & (2k-1)(k-1)\prod_{i=0}^{2k-4}(n-3-i)gc^2_g\dot R_gh+
(2k-1)\times \\
& & \,\,\,\,\,\,\,\,\,\,\,\,\,\,\,\,\,\,\,\,\,\,\,\,\,\,\,\,\,\,\,\,\,\,\,\,\,\times \prod_{i=0}^{2k-3}(n-3-i)c_g\dot R_gh\\
& = & (2k-1)\prod_{i=0}^{2k-4}(n-3-i)\big((k-1)gc^2_g\dot R_gh+\\
& &  \qquad\qquad \qquad\qquad\qquad+(n-2k)c_g\dot R_gh\big),
\end{eqnarray*}
so that
\begin{eqnarray}
c^{2k-1}_gg^{2k-2}\dot
R_gh=\frac{(2k-1)!(n-3)!}{(n-2k)!}\left((k-1)(c^2_g\dot R_gh)g+
(n-2k)c_g\dot R_gh\right),\nonumber\label{aa}
\end{eqnarray}
which gives
\begin{eqnarray}
B_gh=C_{n,k}\mu^{k-1}\left((k-1)(c_g^2\dot R_g h)g+
(n-2k)c_g\dot R_g h\right). \label{bb}
\end{eqnarray}

This shows that, in the constant curvature case, $B_gh$
is completely determined by the first and second contractions of the projection of $\dot R_g h$ on ${\mathcal C}^2(X)$. Now, these contractions have essentially been computed in Lemma \ref{aux}. In fact, inserting the expressions there in (\ref{varl}) we find that

\begin{equation}\label{gg}
c_g\dot R_gh=
\frac12 \left(\nabla^*\nabla h-
\nabla d\mathsf {tr}_gh-2\delta^{*}_g\delta_g h+{\mathcal R}_g\circ h+h\circ {\mathcal R}_g \right)
\end{equation}
and
\begin{equation}
c^2_g \dot R_gh=\Delta_g \mathsf {tr}_gh+\delta_g (\delta_g h)+ \langle{\mathcal R}_g,h\rangle\label{rr},
\end{equation}
so that by (\ref{bb}) we get

\begin{eqnarray}
B_g
h & = &  C_{n,k}\mu^{k-1}\bigg((k-1)\left(\Delta_g \mathsf {tr}_gh+\delta_g (\delta_g h)+ \langle {\mathcal R}_g,h\rangle\right)g +  \nonumber\\
& &\qquad +\frac{n-2k}2 \Big(\nabla^*\nabla h-
\nabla d\mathsf {tr}_gh-\\
& & \qquad\qquad-2\delta^{*}_g\delta_g h+({\mathcal R}_g\circ h+h\circ {\mathcal R}_g) \Big)\bigg).\nonumber
\label{ss}
\end{eqnarray}
Finally, using Remark \ref{newr} to specialize to the constant curvature case, together with (\ref{varl}) and (\ref{zz}) we obtain an expression for $\dot{\mathcal R}^{(2k)}_g$ and hence for ${\mathcal L}_g^{(2k)}=\dot{\mathcal R}^{(2k)}_g-\lambda_k$, namely,

\begin{eqnarray*}
{\mathcal L}_g^{(2k)}h& = & C_{n,k}\mu^{k-1}\bigg( (k-1)\Big(\Delta_g
{\rm tr}_g(h)+\delta_g (\delta_g h)-\\
& & \qquad\qquad\qquad\qquad\qquad -\frac{(k(n-3)+1)}{k-1}\mu {\rm tr}_gh\Big)g+ \\
&&\qquad + \frac{n-2k}2\Big( \nabla^*\nabla h-
\nabla d{\rm tr}_gh-2\delta^{*}_g\delta_g h + 2\mu h \Big)
\bigg).
\end{eqnarray*}

In the sequel we shall use two special cases of this formula, namely,

\begin{eqnarray}
{\mathcal L}_g^{(2k)}\Big|_{{\rm tr}_g^{-1}(0)} & =  & C_{n,k}\mu^{k-1}\bigg( \frac{n-2k}2\nabla^*\nabla
-(n-2k)\delta^{*}_g\delta_g +\nonumber \\
& & \qquad\qquad+(k-1)\delta_g \delta_g (\cdot)g + (n-2k)\mu
\bigg)\label{trace}
\end{eqnarray}
and
\begin{equation}\label{sp2}
{\mathcal L}_g^{(2k)}\Big|_{{\mathcal I}_g}= (n-2k)C_{n,k}\mu^{k-1}\bigg( \frac12\nabla^*\nabla
+ \mu
\bigg).
\end{equation}
In particular, since $\nabla^*\nabla$ is elliptic, (\ref{sp2}) already shows that $\dim\mathfrak{M}_{[g]}^{(2k)}<+\infty$ if $\mu\neq 0$, thus proving the first part of Theorem \ref{main}.

To show that spherical space forms are infinitesimally $2k$-non-deformable, we assume $\mu>0$ and integrate (\ref{sp2}) with $0\neq h\in{\mathcal I}_g$, so that
\begin{eqnarray*}
0  & = & \int_X\langle {\mathcal L}_g^{(2k)}h,h\rangle\nu_g \\
   & = & (n-2k)C_{n,k}\mu^{k-1}\left(\frac{1}{2}\int_X|\nabla h|^2\nu_g+\mu\int_X|h|^2\nu_g\right)\\
   & \geq & (n-2k)C_{n,k}\mu^k\int_X|h|^2\nu_g\\
   & > & 0,
\end{eqnarray*}
a contradiction, thus showing that $\mathfrak{M}_g^{(2k)}= \ker {\mathcal L}_g^{(2k)}|_{{\mathcal I}_g}$ is trivial.

\section{The proof of Theorems \ref{main2} and \ref{main3}}

As remarked in the Introduction, it is not hard to combine the above information on ${\mathcal L}_g^{(2k)}|_{{\mathcal I}_g}$, namely, that it is a self-adjoint positive elliptic operator, with results in \cite{K1} to check that spherical space forms are $2k$-non-deformable indeed. However, we will be able to prove a much stronger result. The idea here is to adapt an argument in \cite{Be}, page 351, which consists in replacing the operator ${\mathcal L}_g^{(2k)}|_{{\rm tr}^{-1}(0)}$, which is obviously non-elliptic, by an elliptic one. In order to implement this strategy, we need some more notation.

First, the fact that the $2k$-Einstein-Lovelock tensor ${\mathcal T}_g^{(2k)}$ is divergence free, for any metric $g$, can be expressed as
\begin{equation}\label{free}
\delta_g{\mathcal R}_g^{(2k)}+(2k-1)!d{\mathcal S}_g^{(2k)}=0.
\end{equation}
We now introduce the functional ${\mathcal G}:{\mathcal M}_1(X)\to {\mathcal C}^1(X)$,
$$
{\mathcal G}(g)={\mathcal R}_g^{(2k)}-\frac{(2k)!}{n}{\mathcal A}^{(2k)}(g)g,
$$
and the $(2k)$-Bianchi operator $\beta^{(2k)}_g:{\mathcal C}^1(X)\to \Omega^1(X)$,
$$
\beta^{(2k)}_g=\delta_g+\frac{1}{2k}d{\mathsf {tr}}_g.
$$
The following lemma is a direct consequence of the definitions.

\begin{lemma}\label{dr}
The following properties hold:
\begin{enumerate}
 \item $g$ is $2k$-Einstein if and only if ${\mathcal G}(g)=0$;
 \item If $g$ is $2k$-Einstein then $\dot{{\mathcal G}}_g={\mathcal L}_g^{(2k)}$ on ${\mathcal C}^1(X)$.
 \item For any $g$, $\beta^{(2k)}_g{{\mathcal G}}(g)=0$. In particular, if $g$ is $2k$-Einstein,
 $\beta^{(2k)}_g{{\mathcal L}}^{(2k)}_g=0$.
\end{enumerate}
\end{lemma}

The identity $\beta^{(2k)}_g{{\mathcal L}}^{(2k)}_g=0$ means that ${{\mathcal L}}^{(2k)}_g$ is not surjective but instead satisfies a first order differential equation coming  from the diffeomorphism invariance of the $2k$-Einstein condition. This is of course a serious obstruction to using the Implicit Function Theorem to probe the local structure of the moduli space. As a way to overcome this we use (\ref{trace}) and introduce, for $h\in {\mathsf {tr}}_g^{-1}(0)$,
the elliptic operator
\begin{eqnarray*}
\tilde{\mathcal L}_g^{(2k)}h & = & {\mathcal L}_g^{(2k)}h+(n-2k)C_{n,k}\mu^{k-1}\delta^*_g\delta_gh-
(k-1)C_{n,k}\mu^{k-1}(\delta_g\delta_gh)g\\
 & = & (n-2k)C_{n,k}\mu^{k-1}\left(\frac{1}{2}\nabla^*\nabla h+\mu h \right).
\end{eqnarray*}
Notice that the identity ${\mathsf {tr}}_g\nabla^*\nabla h=\Delta_g{\mathsf {tr}}_gh$ gives
$$
{\mathsf {tr}}_g\tilde{\mathcal L}_gh=(n-2k)C_{n,k}\mu^{k-1}\left(\frac{1}{2}\Delta_g{\mathsf {tr}}_gh
+\mu{\mathsf {tr}}_gh\right),
$$
which implies
$$
\int_X{\mathsf {tr}}_g\tilde{\mathcal L}_gh\,\nu_g=(n-2k)C_{n,k}\mu^{k}
\int_X{\mathsf {tr}}_gh\,\nu_g,
$$
and recalling that
$$
T_g{\mathcal M}_1(X)=\{h\in{\mathcal C}^1(X);\int_X \mathsf {tr}_gh\,\nu_g=0\},
$$
this shows that $\tilde{\mathcal L}_g(T_g{\mathcal M}_1(X))$ is closed.

We now look at the constraints posed on $\tilde{\mathcal L}_g$ by diffeomorphism invariance. Using Proposition \ref{dr} and the identity ${\mathsf {tr}}_g\delta_g^*\omega=-\delta_g\omega$, $\omega\in\Omega^1(X)$, we compute:
\begin{eqnarray*}
\beta^{(2k)}_g\tilde{\mathcal L}_gh & = & C_{n,k}\mu^{k-1}\left((n-2k)\left(\delta_g
\delta_g^*\delta_gh+\frac{1}{2k}d{\mathsf {tr}}_g\delta_g^*\delta_gh\right.
\right)-\\
&  & \,\,\,\,\,\,\,\,\,\,\,\,\,-(k-1)\left.\left(
\delta_g(\delta_g\delta_gh)g+\frac{1}{2k}d{\mathsf {tr}}_g(\delta_g\delta_gh)g
\right)\right)\\
& = & C_{n,k}\mu^{k-1}\left((n-2k)\left(\delta\delta_g^*-\frac{1}{2k}d\delta\right)(\delta_gh)\right.-\\
& & \,\,\,\,\,\,\,\,\,\,\,\,\,\,\,\left. -(k-1)\left(-d(\delta_g\delta_gh)+\frac{n}{2k}d\delta_g\delta_gh\right) \right)\\
& = &
(n-2k)C_{n,k}\mu^{k-1}\left(\delta_g\delta_g^*-\frac{1}{2}d\delta_g\right)(\delta_gh),
\end{eqnarray*}
so that
\begin{equation}\label{pup}
\beta^{(2k)}_g\tilde{\mathcal L}_gh=(n-2k)C_{n,k}\mu^{k-1}G_g\delta h,
\end{equation}
where
$$
G_g=\delta_g\delta_g^*-\frac12d\delta_g
$$
is elliptic.
Now, (\ref{pup})
first gives $\tilde{\mathcal L}_g(T_g\mathfrak{S}_g)\subset {\ker}\beta^{(2k)}_g$.
Also, if $h=\tilde{\mathcal L}_g k\in \ker\beta^{(2k)}_g$, $k\in T_g{\mathcal M}_1(X)$, then $\delta_gk\in\ker G$, a finite dimensional space, and
this gives
$$
\tilde{\mathcal L}_g(T_g\mathfrak{S}_g)\subset\tilde{\mathcal L}_g(T_g{\mathcal M}_1(X)\cap\ker \beta^{(2k)}_g)\subset
\tilde{\mathcal L}_g\left(T_g{\mathcal M}_1(X)\cap\delta_g^{-1}\ker G\right).
$$
Since $T_g\mathfrak{S}_g$ is closed and has finite codimension
in $T_g{\mathcal M}_1(X)\cap\delta_g^{-1}\ker G$, it is easy to check that $\tilde{\mathcal L}_g(T_g\mathfrak{S}_g)$
is closed in ${\mathcal L}_g\left(T_g{\mathcal M}_1(X)\cap\delta_g^{-1}\ker G\right)$. Thus, the image $\tilde{\mathcal L}_g(T_g{\mathcal M}_1(X))$ is closed in $\tilde{\mathcal L}_g(T_g{\mathcal M}_1(X)\cap\ker \beta^{(2k)}_g)$, which is closed in ${\mathcal C}^1(X)$.

We  conclude that ${\mathcal L}_g=\dot {\mathcal G}_g:T_g\mathfrak{S}_g\to {\mathcal C}^1(X)$, even though not surjective, has closed range.
Thus, if $p$ is the orthogonal projection of ${\mathcal C}^1(X)$ onto ${\mathcal L}_g(T_g\mathfrak{S}_g)$,
the real analytic composite map $p\circ {\mathcal G}:\mathfrak{S}_g\to {\mathcal L}_g(T_g\mathfrak{S}_g)$ is a submersion at $g$. Thus, $(p\circ {\mathcal G})^{-1}(0)$ is a real analytic manifold near $g$, with $\mathfrak{M}_{[g]}^{(2k)}$ as its tangent space at $g$. On this manifold, the map ${\mathcal G}$ is real analytic so that the pre-moduli space ${\mathcal E}_{2k}(X)_g={\mathcal G}^{-1}(0)$ is a real analytic subset. This completes the proof of Theorem \ref{main3}, and hence of Theorem \ref{main2}.

\section{Further comments and questions}

The results in this note give rise to a few basic questions on $2k$-Einstein structures that we would like to briefly discuss here. At a more basic level, it would be highly desirable to find geometric conditions on a $2k$-Einstein metric in order to have the operator $h\mapsto \dot{\mathcal R}^{(2k)}_gh$ elliptic. Since in general the first term in the right hand side of (\ref{sy}) involves no derivatives of $h$, this amounts to understanding the symbol of the linear operator
\begin{equation}\label{ap}
h\mapsto c^{2k-1}_gR^{k-1}_g\dot R_gh
\end{equation}
in more geometric terms. For example, the constant curvature assumption here implies that this symbol is a multiple of the symbol of the Bochner Laplacian, but a more general ellipticity criterium certainly would be of some use. As  is apparent from (\ref{ap}), the difficulty here is that for $k\geq 2$ the $2k$-Einstein condition is fully nonlinear in the second derivatives of the metric, thus implying that the symbol of the linearized operator depends on second order data, namely, the curvature. In any case,  progress in this issue  could be useful in discussing the eventual rigidity of other classes of  $2k$-Einstein manifolds like sufficiently  pinched manifolds, hyperbolic manifolds and certain classes of symmetric spaces.

\bibliographystyle{amsplain}

\end{document}